# Borel-Weil-Bott Theory for Loop Groups

Constantin Teleman
Saint John's College, Cambridge, CB2 1TP, U.K.
email: teleman@dpmms.cam.ac.uk

Loop groups admit a class of representations (those of positive energy) similar to those of semi-simple groups, for which the subgroups of loops that extend holomorphically to a compact Riemann surface play the role of maximal parabolic subgroups. In this note, I outline the corresponding Borel-Weil-Bott theory, as worked out in [T1], [T2], [T3]. The general ideas, in the analytic setting of Conformal Field Theory, are due to Graeme Segal, while my algebraic methods rely on the works [K], [M], [LS], [TUY]. There is now a body of literature ([BL], [F], [KNR]) devoted to the "Borel-Weil" theory (spaces of holomorphic sections); my contribution mainly concerns the "Bott" part (higher cohomology), except for the uniform determination of the "fusion rules", which is a noteworthy application thereof.

Let $G$ be a complex, simple[1], simply connected Lie group, $\Sigma$ a smooth, affine curve, smoothly compactified to $\Sigma^c$, with points at infinity $w_1,\ldots,w_n$. The group $G^\Sigma$ (or $G[\Sigma]$) of algebraic $G$-valued functions on $\Sigma$ is contained in the product $\hat{L}G$ of its formal completions $\hat{L}_iG$ at the $w_i$. These are copies of the group of $G$-valued formal Laurent series. *Highest-weight representations* (HWR's) of $\hat{L}G$ are tensor products of HWR's of the $\hat{L}_iG$. The factors, which shall all be chosen at the same, negative level $(-h)$, are direct sums of the (negative) energy eigenspaces. The projective cocycle of such HWR's splits uniquely over $G^\Sigma$. Finite sums of algebraic duals of HWR's are the *positive energy representations* (PER's); they are formally complete, for a total-energy filtration.

Label $m$ points $z_1,\ldots,z_m$ on $\Sigma$ by irreducible representations (irreps) $\mathsf{V}_k$ of $G$, and let $G^\Sigma$ act on the space $\mathbf{V} := \mathsf{V}_1 \otimes \ldots \otimes \mathsf{V}_m$ by evaluation at the respective points. Consider:
- $X$, the product of the flag varieties $G((z))/G[[z]]$ of the $\hat{L}_iG$;
- $X_\Sigma$, the generalized flag variety $G^\Sigma \backslash \hat{L}G$;
- $\mathfrak{M}$, the moduli stack of $G$-bundles over $\Sigma^c$;
- $\mathcal{L}$, the generator of $\mathrm{Pic}(\mathfrak{M})$ (cf. [LS]);
- $\mathcal{V} := \bigotimes_k \mathcal{V}_k(z_k)$, where $\mathcal{V}_k(z_k)$ is the restriction to $\mathfrak{M} \cong \mathfrak{M} \times \{z_k\}$ of the $\mathsf{V}_k$-bundle associated to the universal $G$-bundle on $\mathfrak{M} \times \Sigma^c$.

$X$ is an ind-scheme, whereas $X_\Sigma$ is a scheme of infinite type. The "uniformization theorem" in [LS] says that the moduli stack $\mathfrak{M}$ is (étale) equivalent to the quotient stack $G^\Sigma \backslash X$. (This combines A. Weil's adèlic description of moduli spaces of vector bundles with the holomorphic double coset construction of [PS].) Using the subgroup $\hat{L}_+G \subset \hat{L}G$ of formal-holomorphic loops, this is also the quotient stack $X_\Sigma/\hat{L}_+G$. The lifting of $\mathcal{V}$ to $X$ is a trivial bundle, but its fiber $\mathbf{V}$ carries the $G^\Sigma$-action described above. Next, $\mathrm{Pic}(G((z))/G[[z]]) \cong \mathbb{Z}$, and the lift to $X$ of $\mathcal{L}$, which I denote by $\mathcal{L}$ as well, is the product of the positive, generating line bundles on the factors of $X$. Also call $\mathcal{L}$ and $\mathcal{V}$ the lifts to $X_\Sigma$ of the same bundles over $\mathfrak{M}$; their domicile will be contextually clear.

**Theorem 4.** *The cohomology groups $H^*\big(X_\Sigma; \mathcal{L}^{\otimes h} \otimes \mathcal{V}\big)$ are finite sums of HWR's of $\hat{L}G$ at level $(-h)$, with multiplicity space $H^*\big(\mathfrak{M}; \mathcal{L}^{\otimes h} \otimes \mathcal{V} \otimes \mathcal{U}^t\big)$ for the HWR with highest-energy spaces $\mathsf{U}_i$.*

Here, $\mathcal{U}$ is the evaluation bundle analogous to $\mathcal{V}$, but with the irreps $\mathsf{U}_i$ attached to the points $w_i$. The result holds even when $\Sigma$ is singular, but in the smooth case I can be more precise.

---

[1] The extension to semi-simple groups is immediate; see [T3], Sec. V, for the case of non-trivial $\pi_1$.

**Theorem 3.** *The cohomology of $\mathcal{L}^{\otimes h} \otimes \mathcal{V}$ over $\mathfrak{M}$ lives in a single degree $\ell$. It vanishes altogether, if $(\lambda_k + \rho) \cdot \alpha \in (h+c) \cdot \mathbb{Z}$, for any of the highest weights $\lambda_k$ of $\mathsf{V}_k$ and roots $\alpha$ of $\mathfrak{g}$; otherwise, its dimension is given by a "Verlinde factorization formula".*

The degree $\ell$ is the sum of the lengths of the affine Weyl transformations taking the $(\lambda_k + \rho)$ to the positive alcove, at level $(h+c)$; $c$ is the dual Coxeter number of $\mathfrak{g}$.

*Proof of Thm. 4:* The key statement is the highest-weight property; see [T3], Sec. VIII for the proof. Assuming that, let $\mathbf{U} := \bigotimes_i \mathsf{U}_i$, on which $\hat{L}_+ G$ acts by evaluation at the $w_i$. There is a "descent spectral sequence" from $X_\Sigma$ to $\mathfrak{M}$, involving the *algebraic group cohomology* of $\hat{L}_+ G$:

$$E_2^{p,q} = H^p_{\hat{L}_+ G}\left(H^q\left(X_\Sigma; \mathcal{L}^{\otimes h} \otimes \mathcal{V}\right) \otimes \mathbf{U}^t\right) \Rightarrow H^*\left(\mathfrak{M}; \mathcal{L}^{\otimes h} \otimes \mathcal{V} \otimes \mathcal{U}^t\right). \tag{1}$$

For an HWR $\check{\mathbf{H}}$, $H^p_{\hat{L}_+ G}\left(\check{\mathbf{H}} \otimes \mathbf{U}^t\right)$ vanishes for $p > 0$, while $H^0$ equals $\mathbb{C}$, rather than zero, precisely when $\check{\mathbf{H}}$ has highest-energy spaces $\mathsf{U}_i$. The sequence (1), then, collapses at $E_2$, and the abutment is the multiplicity space for this same $\check{\mathbf{H}}$. ♠

*Proof of Thm. 3:* By a theorem of Kumar [K] and Mathieu [M], we have $H^0\left(X; \mathcal{L}^{\otimes h}\right) = \mathbf{H}_{0,h}$ (the "vacuum representation" at level $h$), while the higher cohomology vanishes. The collapse at $E_2$ of the descent spectral sequence for the morphism $X \to \mathfrak{M}$ implies ($p = *, q = 0$)

$$H^*\left(\mathfrak{M}; \mathcal{L}^{\otimes h} \otimes \mathcal{V}\right) \cong H^*_{G[\Sigma]}\left(\mathbf{H}_{0,h} \otimes \mathbf{V}\right), \tag{2}$$

the last term being the group cohomology of $G^\Sigma$. The theorem follows from the next result. ♠

**Theorem 2.** *For any PER $\mathbf{H}$ of $\hat{L}G$ at level $h$, the group cohomology $H^*_{G[\Sigma]}(\mathbf{H} \otimes \mathbf{V})$ is concentrated in degree $\ell$, and satisfies the properties listed in Thm. 3.*

*Remark.* This is modeled on the following rewriting of the Borel-Weil-Bott theorem. Let $P \subset G$ be a parabolic subgroup, $\mathsf{E}$ an irrep of $G$, $\mathsf{F}$ one of $P$. The group cohomology $H^*_P(\mathsf{E} \otimes \mathsf{F})$ of $P$, with coefficients in $\mathsf{E} \otimes \mathsf{F}$, is determined as follows. If $\mathcal{F}$ is the sheaf of sections of the algebraic vector bundle $G \times^P \mathsf{F}$ over $G/P$, the descent spectral sequence from $G$ to $G/P$,

$$E_2^{p,q} = H^p_P\left(H^q(G; \mathcal{O}_G \otimes \mathsf{F})\right) \Rightarrow H^*(G/P; \mathcal{F}), \tag{3}$$

collapses at $E_2$, because $G$ is affine. Using the Peter-Weyl theorem for $G$, we get

$$H^*(G/P; \mathcal{F}) \cong \bigoplus_\mathsf{E} \mathsf{E}^t \otimes H^*_P(\mathsf{E} \otimes \mathsf{F}). \tag{4}$$

BWB says that the left-hand side lives in a single degree, where it gives an irrep of $G$; this amounts to the vanishing of $H^*_P(\mathsf{E} \otimes \mathsf{F})$, except possibly in a single degree. Pure-dimensionality of $H^*_{G[\Sigma]}(\mathbf{H} \otimes \mathbf{V})$ is the corresponding loop group statement. Note, however, that this line of argument does not work for loop groups, for which the Peter-Weyl statement fails.

*Proof of Thm. 2:* The *van Est* spectral sequence relates group cohomology to Lie algebra cohomology and to the cohomology of the topological space underlying $G^\Sigma$. It reads

$$E_2^{p,q} = H^p_{G[\Sigma]}(\mathbf{H} \otimes \mathbf{V}) \otimes H^q(G[\Sigma]) \Rightarrow H^*(\mathfrak{g}[\Sigma]; \mathbf{H} \otimes \mathbf{V}). \tag{5}$$

Its collapse at $E_2$, and the single-dimensionality theorem, follow from the next two theorems. ♠

**Theorem 1.** *$G^\Sigma$ is homotopy equivalent to the group $C^\infty(\Sigma; G)$ of smooth maps to G.*

*Proof:* This follows by equating homotopy types in two descriptions of the stack of holomorphic $G$-bundles. Regarding it as the analytic stack underlying the stack of algebraic $G$-bundles gives the homotopy quotient $X/G^\Sigma$, by the uniformization theorem. On the other hand, the Atiyah-Bott construction [AB] realizes it as the quotient stack of smooth $(0,1)$-connections modulo complex gauge transformations. This gives the homotopy type $C^\infty(\Sigma^c; BG)$, which is the homotopy quotient $\Omega G^{\times n}/C^\infty(\Sigma; G)$. Theorem 1 follows because $X \sim \Omega G^{\times n}$, by the natural inclusion (cf. [PS], Ch. 8). ♠

*Remark.* The map $H^*(L\mathfrak{g}; \mathbb{C}) \to H^*(LG)$ is onto ([PS], Ch. 4), and Thm. 1 implies that the same holds for the left edge-homomorphism in (5). The sequence, therefore, collapses at $E_2$.

**Theorem 0.** [T2, Thm. 2.5] *The cohomology of the Lie algebra $\mathfrak{g}[\Sigma]$ of $\mathfrak{g}$-valued algebraic functions on $\Sigma$, with coefficients in the representation $\mathbf{H} \otimes \mathbf{V}$, is given by*

$$H^*(\mathfrak{g}[\Sigma]; \mathbf{H} \otimes \mathbf{V}) \cong H^\ell(\mathfrak{g}[\Sigma]; \mathbf{H} \otimes \mathbf{V}) \otimes H^{*-\ell}(C^\infty(\Sigma; G)) \tag{6}$$

*The first factor on the right satisfies the properties listed in Thm. 3.*

*Remark.* The second factor is $H^*(G) \otimes H^*(\Omega G^{\times N})$, where $N = 2g + n - 1$.

*Proof:* The proof goes by induction over the genus. For the inductive step, let $\Sigma$ degenerate to a curve $\Sigma_0$ of genus $(g-1)$, with one node and normalization $\tilde{\Sigma}_0$. Shapiro's lemma gives

$$H^*_{G[\Sigma_0]}(\mathbf{H} \otimes \mathbf{V}) \cong H^*_{G[\tilde{\Sigma}_0]}\left(\mathrm{Ind}^{G[\tilde{\Sigma}_0]}_{G[\Sigma_0]}(\mathbf{H} \otimes \mathbf{V})\right) \cong \bigoplus_{\mathsf{U}} H^*_{G[\tilde{\Sigma}_0]}\left(\mathbf{H} \otimes \mathbf{V} \otimes \mathsf{U}(x') \otimes \mathsf{U}^t(x'')\right), \tag{7}$$

using evaluation representations at $x'$ and $x''$, the liftings in $\tilde{\Sigma}_0$ of the node. By inductive assumption in Thm. 2, $H^\ell_{G[\Sigma_0]}(\mathbf{H} \otimes \mathbf{V})$ is given by Verlinde factorization in genus $(g-1)$, and one shows as in [T2], Prop. 3.8, that the entire group cohomology is $H^\ell_{G[\Sigma_0]}(\mathbf{H} \otimes \mathbf{V}) \otimes H^*(\Omega G)$. (This uses Bott's basis for $H_*(\Omega G)$, indexed by alcoves in the positive Weyl chamber.) Further,

$$H^*(\mathfrak{g}[\Sigma_0]; \mathbf{H} \otimes \mathbf{V}) \cong H^*_{G[\Sigma_0]}(\mathbf{H} \otimes \mathbf{V}) \otimes H^*(G[\Sigma_0]) \cong H^\ell_{G[\Sigma_0]}(\mathbf{H} \otimes \mathbf{V}) \otimes H^{*-\ell}(G \times \Omega G^{\times N}), \tag{8}$$

by the van Est spectral sequence, the last isomorphism following from $G[\Sigma_0] \sim G \times \Omega G^{\times(N-1)}$. Finally, consider the *relative* cohomology $H^*(\mathfrak{g}[\Sigma_0], \mathfrak{g}; \mathbf{H} \otimes \mathbf{V})$: it loses the factor $H^*(G)$, so it lives only in degrees of the same parity; this implies that relative cohomology (and then, the absolute one) is rigid under the degeneration of $\Sigma$ to $\Sigma_0$, and is given as in Thm. 0.

In genus zero, further degeneration reduces us to the case of a single puncture, treated in [T1], Thm. 0. The Lie algebra cohomology[2] $H^*(\mathfrak{g}[z^{-1}], \mathfrak{g}; \mathbf{H} \otimes \mathbf{V})$ is resolved by a Koszul complex $\mathcal{C}^*$, sitting inside the Dolbeault complex for $\mathcal{L}^{\otimes h} \otimes \mathcal{V}$ over the flag variety $X_0 := G[z^{-1}]\backslash G((z))$. $\mathcal{C}^*$ carries a natural Hilbert space topology, involving the standard hermitian norms on $\mathbf{H}$ and $\mathbf{V}$, and the Kähler metric on $X_0$. (Both require us to fix the unit disk in the affine line). The Weitzenböck formula (also known here as Nakano's identity) expresses the $\bar{\partial}$-Laplacian as a non-negative operator (the $(1,0)$-Laplacian) plus a zero-order term, related to the hermitian curvature of $\mathcal{L}^{\otimes h} \otimes \mathcal{V}$ and to the Ricci curvature of $X_0$. All of these preserve $\mathcal{C}^*$.

When the shifted highest weights $(\lambda_k + \rho)$ lie in the positive alcove, and when the $z_k$ are just outside the unit disk, but not near each other, it turns out that the curvature term is positive on positive-degree forms. There is, then, no higher Lie algebra $L^2$-cohomology[3]. A direct computation ([T1], Sec. 3.3) shows that any cohomology class in the algebraic Koszul complex — the direct product of energy eigenspaces — has a Hilbert space representative; but this must lie in the range

---

[2] This differs from the convention in [T1], where $\mathfrak{g}[z]$ was used instead.

[3] Under the stated assumptions, Nakano's identity also ensures that the Hilbert space range of $\bar{\partial}$ is closed; the argument in Sec. 3.3 of [T1] was excessively cautious.

of $\bar\partial$, so all higher "honest" Lie algebra cohomology vanishes as well. In the algebraic context, the restriction on the placement of the $z_k$ can be readily lifted. This yields the case $\ell = 0$ of Thm. 0; the general case follows from the usual reflection argument of Bott ([T1], Sec. 3.5). ♠

Finally, recall the connection with the fusion rules, originally conjectured by Segal. Bott defines a linear *holomorphic induction* map $I_h$, from the representation ring of $G$ to the free Abelian group on HWR's of the loop group at level $(-h)$: V is sent to the Euler characteristic of $\mathcal{L}^{\otimes h} \otimes \mathcal{V}$ over $X_0$ ($m = 1$, $V_1 := V$, $\Sigma$ is the affine line). The choice of a marked point $z_1$ is irrelevant; when V is irreducible, $I_h(V)$ (if non-zero) is $\pm$ a HWR of $LG$, with virtual character given by the Kac formula, applied to the highest weight of V.

Given $m$ HWR's of of $LG$, with highest-energy spaces the $G$-irreps $V_1,\ldots,V_m$, Segal defines their *fusion product* as the space of global sections of $\mathcal{L}^{\otimes h} \otimes \mathcal{V}$ over $X_0$. This depends on the choice of $m$ distinct points on the line, but its HWR-multiplicities do not, and can be equated with other definitions of the fusion coefficients, using the degree-zero case of Thms. 2, 3 and 4 (cf. also [T2], Prop. 2.10). (The multiplicity *spaces*, of course, carry extra structure, namely a flat projective connection and, conjecturally, a distinguished, compatible hermitian metric.) The following was, I think, first conjectured by Bott.

**Corollary 5.** *Holomorphic induction is a ring homomorphism, taking tensor product to fusion.*

*Proof:* By Thms. 3 and 4, the fusion product is really the Euler characteristic of $\mathcal{L}^{\otimes h} \otimes \mathcal{V}$ over $X_0$. This does not change when $\mathcal{V}$ is deformed by moving the marked points $z_1,\ldots,z_m$ to $\infty$, so it equals the desired $I_h(V_1 \otimes \ldots \otimes V_m)$. To see this rigidity of the characteristic, note that the fusion multiplicities are given by Lie algebra cohomology groups (Thm. 0); but ([T1], Cor. 3.2.7) these can be resolved by a finite-dimensional complex, where only the differentials depend on the $z_k$. ♠

*Remark.* This proof (dating from late '93, in Lie algebra language) realizes an older, heuristic argument of Segal's, and was probably the first uniform proof published ([T1], Thm.1); although several other arguments, partial or exhaustive (e.g. [F], [Fi]), are slightly prior to it. There seems to be no simple uniform proof: [Fi] uses results of Gelfand-Kazhdan and Kazhdan-Lusztig.